\documentclass[12pt]{article}
\usepackage{amsfonts,latexsym,amssymb}
 \newcommand{\beqn}{\begin{eqnarray}}
 \newcommand{\eeqn}{\end{eqnarray}}
 \newcommand{\be}{\begin{equation}}
 \newcommand{\ee}{\end{equation}}
 \newcommand{\ba}{\begin{array}}
 \newcommand{\ea}{\end{array}}
 \newcommand{\pa}{\partial}
 
 \newcommand{\ci}{\cite}
 \newcommand{\la}{\label}

\newcommand{\Om}{\Omega}

\newcommand{\ga}{\gamma}

\newcommand{\na}{\nabla}
\newcommand{\om}{\omega}

 \newcommand{\de}{\delta}

 \newcommand{\De}{\Delta}

\def\R{{\rm I\kern-.1567em R}}
\def\M{{\rm I\kern-.1567em M}}
\def\intp{{\rm \int{\hspace{-4mm}}-}}
\def\tr{{\rm tr}}

\def\div {{\rm div}}

\def\dist{{\rm dist}}
\def\ess{{\rm ess}}

 \newtheorem{theorem}{Theorem}[section]

 \newtheorem{lemma}[theorem]{Lemma}
 
 \newtheorem{remark}[theorem]{Remark}

 \newtheorem{pro}[theorem]{Proposition}

\begin{document}

\begin{center} {\Large\bf   On Smoothness of $L_{3,\infty}$-Solutions to the
Navier-Stokes Equations up to  Boundary }\\
  \vspace{1cm}
 {\large
  G. Seregin}

 \end{center}

 \vspace{1cm}
 \noindent
 {\bf Abstract } We show that $L_{3,\infty}$-solutions to
  the three-dimensional
Navier-Stokes equations near a flat part of the boundary are
smooth.

 \vspace {1cm}

\noindent {\bf 1991 Mathematical subject classification (Amer.
Math. Soc.)}: 35K, 76D.

\noindent
 {\bf Key Words}: boundary regularity, the Navier-Stokes
 equations, suitable weak solutions, backward uniqueness.

\setcounter{equation}{0}
\section{Introduction  }

In the present paper, we are going to prove smoothness of the
so-called $L_{3,\infty}$-solutions to the Navier-Stokes equations
up to a flat part of the boundary. In particular,  Theorem
\ref{it1} proved below implies the result  announced in \ci{ESS3}.
It was stated there that $L_{3,\infty}$-solutions to the initial
boundary value problems for the Navier-Stokes equations in a half
space are smooth if the initial data are smooth. As in the case of
the Cauchy problem, we deduce this statement from the theorem on
local regularity of $L_{3,\infty}$-solutions near a flat part of
the boundary.

The main idea how to treat boundary regularity of
$L_{3,\infty}$-solutions is similar
to the case of interior regularity: reduction to  a backward
uniqueness problem for the heat operator, see \ci{SS1}, \ci{ESS2},
and \ci{ESS4}. The second part of such analysis has been already
done in \ci{ESS3}, where the backward uniqueness result for the
heat operator in a half space was established.

However,  serious difficulties occur if we scale and blow up the
Navier-Stokes equations
at singular boundary points. In particular, since
$L_{3,\infty}$-norm is invariant with respect to the natural
scaling,
the global $L_{3,\infty}$-norm of the blow-up velocity is bounded.
In the interior case, we were able to prove global boundedness of
$L_{\frac 32,\infty}$-norm of the the blow-up pressure. We do not
know whether the same is true near the boundary. If it would be
so,  the proof of boundary regularity could be essentially
simplified.
 Unfortunately, we cannot
even show that there is a reasonable global norm of the blow-up
pressure which is finite. This makes our proof a  bit tricky. Key
points are Lemma \ref{4l1} and suitable decomposition of the
pressure.

The main result of the paper is as follows.
\begin{theorem}\la{it1} Let a pair of functions $v$ and $p$ has the following
differentiability properties:
\begin{eqnarray}\la{i1}
  v\in L_{2,\infty}(Q^+)\cap W^{1,0}_{2}(Q^+)\cap
   W^{2,1}_{ \frac 98,\frac 32}(Q^+),\qquad p\in
    W^{1,0}_{\frac 98,\frac 32}(Q^+).
\end{eqnarray}
Here, $Q^+=\{z=(x,t)\,\,\|\,\,|x|<1,\,\,x_3>0,\,\,-1<t<0\}$.

Suppose  that $v$ and $p$ satisfy the Navier-Stokes equations a.e.
in $Q^+$, i.e.: \be\la{i2}\left.\begin{array}{c}
\pa_tv+\div\,v\otimes v-\De\,v=-\na\,p
    \\
  \div\,v=0
\end{array}\right\}\quad\mbox{in}\quad Q^+\ee
and the boundary condition
\begin{equation}\label{i3}
    v(x,t)=0,\qquad x_3=0\,\,\mbox{and}\,\,-1\leq t\leq 0.
\end{equation}
Assume, in addition, that
\begin{equation}\label{i4}
    v\in L_{3,\infty}(Q^+).
\end{equation}
Then $v$ is H\"older continuous in the closure of the set
$$Q^+(1/2)=\{z=(x,t)\,\,\|\,\,|x|<1/2,\,\,x_3>0,\,\,-(1/2)^2<t<0\}.$$
\end{theorem}
Explanations why conditions (\ref{i1}) are natural can be found in
 paper \ci{S3}, see Theorem 2.2  there. We just briefly note that any
weak Leray-Hopf solution  to  initial boundary value problems in a
half space together with the associated pressure satisfies
(\ref{i1}). So, the real additional assumption of Theorem
\ref{it1} is condition (\ref{i4}).

\setcounter{equation}{0}
\section{Some preliminary estimates for the pressure }

We denote by $\mathbb M^3$ the space of all real $3\times 3 $
matrices.  Adopting summation over repeated Latin indices, running
from 1 to 3,  we shall use the following notation
$$
u\cdot v =u_iv_i, \quad | u|=\sqrt{u\cdot u}, \quad u=(u_i)\in
\mathbb R^3, \,  v=(v_i)\in \mathbb R^3;
$$
$$
A:B=\tr A^* B=A_{ij}B_{ij}, \quad |A|=\sqrt{A:A},
$$
$$
A^*=(A_{ji}),\quad \tr A= A_{ii}, \quad A=(A_{ij})\in\mathbb M^3,
\,B=(B_{ij})\in\mathbb M^3;
$$
$$
u\otimes v=(u_iv_j)\in\mathbb M^3, \quad Au=(A_{ij}u_j)\in\mathbb
R^3, \quad u,v\in\mathbb R^3, \,A\in \mathbb M^3.
$$

Let $\om$ be a domain in some finite-dimensional space. We denote
by $L_m(\om)$ and $W^l_m(\om)$ the known Lebesgue and Sobolev
spaces. The norm of the space $L_m(\om)$ is denoted by
$\|\cdot\|_{m,\om}$. If $m=2$, then we use the abbreviation
$\|\cdot\|_{\om}\equiv\|\cdot\|_{2,\om}$.

Let $T$ be  a positive parameter, $\Om$ be a domain in ${\mathbb
R}^3$. We denote by $Q_T\equiv\Om\times]0,T[$ the space-time
cylinder. Space-time  points are denoted by $z=(x,t)$,
$z_0=(x_0,t_0)$, and etc.

For summable in $Q_T$ scalar-valued, vector-valued, and
tensor-valued functions, we shall use the following differential
operators $$ \pa_t v =\frac {\pa v}{\pa t}, \quad v_{,i}=\frac
{\pa v}{\pa x_i}, \quad \nabla p=(p_{,i}), \quad \na \,u =
(u_{i,j}), $$ $$ \div \,v =v_{i,i}, \quad\div \,\tau
=(\tau_{ij,j}),\quad \De \,u= \div\,\na\,u, $$ which are
understood in the sense of distributions. Here, $x_i,\,i=1,2,3$,
are the Cartesian coordinates of a point $x\in\mathbb R^3$, and
$t\in ]0,T[$ is a moment of time.

Let $L_{m,n}(Q_T)$ be the space of measurable $\mathbb R^l$-valued
functions with the following norm
$$
  \|f\|_{m,n,Q_T} =
\left\{  \begin{array}c
\big(\int\limits_0^T\|f(\cdot,t)\|^n_{m,\Om}\,dt\big)^{\frac
1n}\,,
 \qquad n \in
[1,+\infty [
\\
\mbox{ess} \sup_{t\in [0,T]}
 \|f(\cdot,t)\|_{m,\Om}\,  \,, \qquad n =+\infty .
  \end{array}
\right.
$$
Now, we can define the following Sobolev spaces with the mixed
norm:
$$W^{1,0}_{m,n}=\{v\in L_{m,n}(Q_T)\,\,\|\,\, \|v\|_{m,n,Q_T}+
\|\na\,v\|_{m,n,Q_T}<+\infty\},$$
$$W^{2,1}_{m,n}=\{v\in L_{m,n}(Q_T)\,\,\|\,\, \|v\|_{m,n,Q_T}+
\|\na\,v\|_{m,n,Q_T}+\|\na^2 v\|_{m,n,Q_T}$$ $$+\|\pa_t
v\|_{m,n,Q_T}<+\infty\}.$$

Setting $x^\prime=(x_1,x_2)\in \mathbb R^2$, we
 introduce the additional notation:
$$
B(x_0,R)\equiv \{x\in \mathbb R^3 \,\,\|\,\, |x-x_0| < R\},
$$
$$
B^+(x_0,R)\equiv \{x\in  B(x_0,R)\,\,\|\,\, x=(x^\prime,x_3),
\quad x_3 > x_{03}\},
$$
$$
B(\theta)\equiv B(0, \theta),\quad B\equiv B(1), \quad
B^+(\theta)\equiv B^+(0, \theta),\quad B^+\equiv B^+(1),
$$
$$\Gamma(x_0,R)\equiv \{x\in B(x_0,R)\,\,\|\,\,x_3=x_{30}\},\quad
\Gamma(\theta)\equiv\Gamma(0,\theta),\quad\Gamma\equiv \Gamma(1),
$$
$$
Q(z_0,R)\equiv B(x_0,R)\times ]t_0-R^2,t_0[, \quad z_0=(x_0,t_0),
$$
$$
Q^+(z_0,R)\equiv B^+(x_0,R)\times ]t_0-R^2,t_0[,
$$
$$
Q(\theta)\equiv Q(0,\theta), \quad Q\equiv Q(1),\quad
Q^+(\theta)\equiv Q^+(0,\theta), \quad Q^+\equiv Q^+(1).
$$

Various mean values of integrable functions are denoted as follows
 $$ [p]_{\Om}(t)\equiv
\intp_{\Om} p(x,t)\,dx \equiv\frac 1{|\Om|}\int\limits_{\Om}
p(x,t)\,dx, $$ $$ (v)_{\om}\equiv\intp_{\om} v\,dz\equiv\frac
1{|\om|}\int \limits_{\om} v\,dz. $$

We denotes by $c$ all universal positive constants.

In this section, we shall prove a couple of propositions about the
pressure in the Navier-Stokes equations provided that conditions
(\ref{i4}) holds. To this end, we are going to use two results
established in \ci{S2} and \ci{S3}. For the reader convenience,
they are  formulated below. Moreover,  since the first lemma is
slightly different from Lemma 3.1 in \ci{S2}, we shall prove it.

\begin{lemma}\la {2l1}Let $v\in L_3(Q(z_0,R))$ and
$p\in L_\frac 32(Q(z_0,R))$ satisfy the Navier-Stokes equations in
the sense of distributions. Then, for $0<r\leq\rho\leq R$, we have
\begin{equation}\label{21}
    D(z_0,r;p)\leq c\Big[\Big(\frac r\rho\Big)^\frac 52D(z_0,\rho;p)
+\Big(\frac \rho r\Big)^2  C(z_0,\rho;v)  \Big],
\end{equation} where
$$C(z_0,r;v)\equiv \frac
1{r^2}\int\limits_{Q(z_0,r)}|v|^3\,dz,\qquad D(z_0,r;p) \equiv
\frac 1{r^2}\int\limits_{Q(z_0,r)}|p-[p]_{B(x_0,r)}|^\frac 32\,dz.
$$
\end{lemma}
\begin{lemma}\la {2l2}
Let a pair of functions $v$ and $p$ satisfy the following
conditions. They have the differentiability properties
\begin{eqnarray}\la{22}
 \nonumber v\in L_{2,\infty}(Q^+(z_0,R))\cap W^{1,0}_{2}(Q^+(z_0,R))\cap
   W^{2,1}_{ \frac 98,\frac 32}(Q^+(z_0,R)),\\
   \qquad p\in
    W^{1,0}_{\frac 98,\frac 32}(Q^+(z_0,R)).
\end{eqnarray}
The pair $v$ and $p$ satisfies the Navier-Stokes equations a.e. in
$Q^+(z_0,R)$ and the boundary condition
\begin{equation}\label{23}
    v(x,t)=0,\qquad x_3=x_{30}\,\,\mbox{and}\,\,\,t_0-R^2\leq t\leq t_0.
\end{equation}
For a.a. $t\in ]t_0-R^2, t_0[$ and for all nonnegative cut-off
functions $\varphi\in C_0^\infty(\mathbb R^4)$, vanishing in a
neighborhood of the parabolic boundary
$$\pa^\prime Q(z_0,R)=B(x_0,R)\times\{t=t_0-R^2\}\cup \pa B(x_0,R)\times
[t_0-R^2, t_0]$$ of the cylinder $Q(z_0,R)$, $v$ and $p$ satisfy
the local energy inequality
\begin{eqnarray}\la{24}
 \nonumber \int\limits_{B(x_0,R)}\varphi(x,t)|v(x,t)|^2\,dx+2\int
  \limits_{B(x_0,R)\times ] t_0-R^2,t[}\varphi |\na\,v|^2\,dxdt^\prime
    \\
  \leq \int
  \limits_{B(x_0,R)\times ] t_0-R^2,t[}\Big[|v|^2(\pa_t\varphi+\na\varphi)
  +v\cdot\na\varphi(|v|^2+2p)\Big]\, dxdt^\prime.\end{eqnarray}
  Then, for any $0<r\leq\rho\leq R$, we have
  \begin{eqnarray}\la{25}
    \nonumber D^+_1(z_0,r;p) \leq c\Big\{\Big(\frac r\rho\Big)^2
    \Big[D^+_1(z_0,\rho;p)+(E^+)^\frac 34(z_0,\rho;v)\Big]\\
    + \Big(\frac \rho r\Big)^\frac 32(A^+)^\frac 12(z_0,\rho;v)
    E^+(z_0,\rho;v)\Big\},
  \end{eqnarray}
  where
  $$A^+(z_0,\rho;v)\equiv \ess \sup\limits_{t_0-R^2<t<t_0}\frac 1\rho
  \int\limits_{B^+(x_0,\rho)}|v(x,t)|^2\,dx,$$
  $$E^+(z_0,\rho;v)\equiv\frac 1\rho
  \int\limits_{Q^+(z_0,\rho)}|\na\,v|^2\,dz,$$
  $$D^+_1(z_0,\rho;p)\equiv \frac
  1{\rho^2}\int\limits_{t_0-\rho^2}^{t_0}\Big( \int\limits_{B^+(x_0,\rho)}
  |\na\,p|^\frac 98\,dx\Big)^\frac 43\,dt.$$
  \end{lemma}
\begin{remark}\la{2r3}Lemma \ref{2l2} was proved in \ci{S3}, see Lemma
7.2 there.\end{remark}
\begin{remark}\la{2r4}According to the definition introduced in
\ci{S3}, see Definition 2.1 there, the pair $v$ and $p$,
satisfying condition (\ref{22})--(\ref{24}), is called a suitable
weak solution to the Navier-Stokes equations in $Q^+(z_0,R)$ near
$\Gamma(x_0,R)\times [t_0-R^2,t_0]$.\end{remark}

\textsc{Proof of Lemma \ref{2l1}} We just outline our proof
because it is essentially the same as the proof of Lemma 3.1 in
\ci{S2}.

For  a.a. $t\in ]t_0-\rho^2,t_0[$, the pressure $p$  meets the
equation
$$\De p(\cdot,t)=-\div\,\div\,v(\cdot,t)\otimes
v(\cdot,t)\qquad\mbox{in}\quad B(x_0,\rho)$$ in the sense of
distributions. We decompose it so that
$$p=p_1+p_2,$$
where $p_1$ is defined as follows:
$$\int\limits_{B(x_0,\rho)}p_1(x,t)\De\varphi(x)\,dx=
\int\limits_{B(x_0,\rho)}v(x,t)\otimes v(x,t):\na^2(x)\,dx$$ for
any $\varphi\in W^2_3(\mathbb R^3)$ such that $\varphi=0$ on $\pa
B(x_0,\rho)$. Regarding $p_2$, we have
\begin{equation}\label{26}
    \De p(\cdot,t)=0\qquad\mbox{in}\quad B(x_0,\rho)
\end{equation}
for  a.a. $t\in ]t_0-\rho^2,t_0[$. By the known regularity theory
results,
\begin{equation}\label{27}
\int\limits_{B(x_0,\rho)}|p_1(x,t)|^\frac 32\,dx\leq
c\int\limits_{B(x_0,\rho)}|v(x,t)|^3\,dx.
\end{equation}

Let $0<r\leq \rho/2$. We have
\begin{eqnarray}\la{28}
\nonumber D(z_0,r;p)\leq c\Big[  D(z_0,r;p_1)+ D(z_0,r;p_2)\Big]  \\
 \nonumber\leq c \Big[\frac  1{r^2}\int\limits_{Q(z_0,\rho)}
 |p_1|^\frac 32\,dz+ D(z_0,r ;p_2)\Big]  \\
  \leq c\Big[\Big(\frac\rho r\Big)^2C(z_0,\rho;v) + D(z_0,r
  ;p_2)\Big].
\end{eqnarray}
Since $p_2$ is a harmonic function, we see that
\begin{eqnarray*}
  \sup\limits_{x\in B(x_0,r)}|p_2(x,t)-[p_2]_{B(x_0,r)}(t)|\leq
  cr\sup\limits_{x\in B(x_0,\rho/2)} |\na\,p_2(x,t)| \\
  \leq cr\frac 1{\rho^4}\int\limits_{B(x_0,\rho)}
  |p_2(x,t)-[p_2]_{B(x_0,\rho)}(t)|dx \\
   \leq c\Big(\frac r\rho\Big)\frac 1{\rho^2}\Big(\int\limits_{B(x_0,\rho)}
  |p_2(x,t)-[p_2]_{B(x_0,\rho)}(t)|^\frac 32\,dx\Big)^\frac 23.
\end{eqnarray*}
Therefore,
\begin{eqnarray}\la{29}
  \nonumber D(z_0,r;p_2)\leq c\frac 1{r^2}\Big(\frac r\rho\Big)^3
  \Big(\frac r\rho\Big)^\frac 32\int\limits^{t_0}_{t_0-r^2}\,dt
  \int\limits_{B(x_0,\rho)}|p_2-[p_2]_{B(x_0,\rho)}|^\frac 32\,dx \\
  \nonumber \leq c  \Big(\frac r\rho\Big)^{3+\frac 32-2}D(z_0,\rho;p_2)  \\
 \nonumber  \leq c\Big(\frac r\rho\Big)^\frac 52\Big[D(z_0,\rho;p)
 +D(z_0,\rho;p_1)\Big]  \\
   \leq c\Big(\frac r\rho\Big)^\frac 52\Big[D(z_0,\rho;p)
 +C(z_0,\rho;v)\Big].
\end{eqnarray}
Combining (\ref{28}) and (\ref{29}), we easily arrived at
(\ref{21}). Lemma \ref{2l1} is proved.

Now, our goal is to prove two auxiliary propositions.
\begin{pro}\la{2p5} Assume that all conditions of Lemma \ref{2l1}
are fulfilled. And let, in addition,
\begin{equation}\label{210}
    \|v\|_{3,\infty,Q(z_0,R)}\leq L<+\infty.
\end{equation}
Then, for any $\ga\in ]0,1[$, there exists a constant $c_1$
depending  on $\ga$ and $L$ only such that, for $0<r\leq R$, we
have
\begin{equation}\label{211}
D(z_0,r;p)\leq c_1\Big[\Big(\frac rR\Big)^{\frac
52\ga}D(z_0,R;p)+1\Big].
\end{equation}
\end{pro}
\textsc{Proof} It can be derived from  (\ref{21})  that:
\begin{equation}\label{212}
D(z_0,\tau^{k+1}R;p)\leq c\Big[\tau^{\frac
52}D(z_0,\tau^kR;p)+\frac 1{\tau^2}L^3\Big]
\end{equation}
for any $0<\tau<1$. We may choose $\tau\in ]0,1[$ so that
$$c\tau^{\frac 52(\ga-1)}\leq 1.$$
Then we find (\ref{211}) from (\ref{212}) just by iterations.
Proposition \ref{2p5} is proved.
\begin{pro}\la{2p6} Assume that all conditions of Lemma \ref{2l2}
are fulfilled. And let, in addition,
\begin{equation}\label{213}
    \|v\|_{3,\infty,Q^+(z_0,R)}\leq L<+\infty.
\end{equation}
Then, for any $\ga\in ]0,1[$, there exists a constant $c_2$
depending  on $\ga$ and $L$ only such that, for $0<r\leq R$, we
have
\begin{equation}\label{214}
D^+_1(z_0,r;p)\leq c_2\Big[\Big(\frac rR\Big)^{
2\ga}D^+_1(z_0,R;p)+1\Big].
\end{equation}
\end{pro}
\textsc{Proof} Let $\rho\leq R/2$. Then local energy inequality
(\ref{24}) gives us the following estimate
\begin{eqnarray*}
  E^+(z_0,\rho;v)\leq c\Big[(C^+)^\frac 23(z_0,2\rho;v)+
(D^+)^\frac 23(z_0,2\rho;p) (C^+)^\frac
13(z_0,2\rho;v)   \\
+C^+(z_0,2\rho;v)\Big],
\end{eqnarray*}
where $$ D^+(z_0,r;p)\equiv\frac 1{r^2}\int\limits_{Q^+(z_0,r)}
|p-[p]_{B^+(x_0,r)}|^\frac 32\,dz.$$ Now, using condition
(\ref{213}) and the embedding theorem, we show
$$E^+(z_0,\rho;v)\leq c\Big[L^2+(D^+_1)^\frac
23(z_0,2\rho;p)L+L^3\Big].$$ And thus, from Lemma \ref{2l2}, see
(\ref{25}), we find
\begin{eqnarray*}
  D^+_1(z_0,r;p)\leq c_3(L)\Big[\Big(\frac r\rho\Big)^2\Big( D^+_1(z_0,2\rho;p)
  +1\Big) \\
  +\Big(\frac\rho r\Big)^3\Big( (D^+_1)^\frac 23(z_0,2\rho;p)+1\Big)\Big] \\
   \leq c^\prime_3(L)\Big[\Big(\frac r\rho\Big)^2 D^+_1(z_0,2\rho;p)+
   \Big(\frac\rho r\Big)^{13}\Big]
\end{eqnarray*}
for any $0<r\leq\rho\leq R/2$. But the latter immediately implies
$$ D^+_1(z_0,r;p)\leq c_4(L)\Big[\Big(\frac r\rho\Big)^2 D^+_1(z_0,\rho;p)+
   \Big(\frac\rho r\Big)^{13}\Big]$$
   for all $0<r\leq\rho\leq R$. Using the same arguments as in the
   proof of the previous proposition, we establish (\ref{214}).
   Proposition \ref{2p6} is proved.

\setcounter{equation}{0}
\section{$\varepsilon$-regularity results for suitable weak solutions}

First, let us show that the pair $v$ and $p$ from Theorem
\ref{it1} forms the so-called suitable weak solution to the
Navier-Stokes equations in $Q^+$ near $\Gamma\times [-1,0]$. The
corresponding definition
was introduced in \ci{S3}. It is a natural modification of the
known definitions
of suitable weak solutions,
discussed in \ci{Sc1}, \ci{CKN}, and \ci{Li}, for the interior
case. In our case, this means that the pair $v$ and $p$ must
subject to the following conditions:
\begin{equation}\label{31}
    v\in L_{2,\infty}(Q^+)\cap W^{1,0}_{2}(Q^+)\cap
   W^{2,1}_{ \frac 98,\frac 32}(Q^+),\qquad p\in
    W^{1,0}_{\frac 98,\frac 32}(Q^+);\end{equation}
    \begin{eqnarray}\la{32}
      v\,\,\mbox{and}\,\,p\,\,\mbox{satisfy the Navier-Stokes
      equations a.e. in}\,\,Q^+;
    \end{eqnarray}
\begin{equation}\label{33}
    v=0\qquad\mbox{on}\quad \Gamma\times [-1,0];
\end{equation}

 for a.a. $t\in ]-1,0[$ and for all nonnegative
 functions  $\varphi \in C^\infty_0(\mathbb R^4)$,
 vanishing a neighborhood of the parabolic boundary
 $\pa' Q$ of $Q$,
   $v$ and $p$ satisfy the inequality
  \begin{eqnarray}\la{34} \nonumber \int\limits_{B^+}\varphi(x,t)|v(x,t)|^2\,dx+2\int
   \limits_{B^+\times ]-1,t[}\varphi |\na\,v|^2\,dxdt^\prime \\
   \leq\int
   \limits_{B^+\times ]-1,t[}\Big[|v|^2(\pa_t\varphi+\De\varphi)+
   v\cdot\na\varphi(|v|^2+2p)\Big]\,dxdt^\prime.
\end{eqnarray}

(\ref{31})--(\ref{33}) hold by the assumptions of Theorem
\ref{it1}. We should just verify that they satisfy local energy
inequality (\ref{34}). To this end, it is sufficient to show that
\begin{equation}\label{35}
    v\in W^{2,1}_\frac 43(Q^+(\tau)),\qquad p\in W^{1,0}_\frac 43(Q^+(\tau))
\end{equation}
for any $\tau\in ]0,1[$. If (\ref{35}) is proved, then (\ref{34})
holds as identity.

 Fix  a  domain $\widetilde{B}$ with smooth boundary such that
$$B^+((1+\tau)/2)\subset\widetilde{B}\subset B^+$$
and consider the following initial boundary value problem
\be\la{36}\left.\begin{array}{c} \pa_tv^1
-\De\,v^1=\widetilde{f}-\na\,p^1
    \\
  \div\,v^1=0
\end{array}\right\}\quad\mbox{in}\quad \widetilde{Q}=
\widetilde{B}\times ]-1,0[\ee
\begin{equation}\label{37}
    v^1|_{\pa^\prime\widetilde{Q}}=0,
\end{equation}
where $\widetilde{f}=-\div\,v\otimes v=-v_iv_{,i}$.  It is easy to
check that
\begin{equation}\label{38}
    \widetilde{f}\in L_\frac 43(Q^+)\cap L_{\frac 98,\frac
    32}(Q^+).
\end{equation}
By the coercive estimates for solutions to the Stokes system, see
\ci{GS} and \ci{MS},
\begin{equation}\label{39}
v^1\in W^{2,1}_{\frac 43}(\widetilde{Q})\cap  W^{2,1}_{\frac
98,\frac 32}(\widetilde{Q}),\qquad p^1 \in   W^{1,0}_{\frac
43}(\widetilde{Q})\cap  W^{1,0}_{\frac 98,\frac
32}(\widetilde{Q}).\end{equation} On the other hand, functions
$v^2=v-v^1$ and $p^2=p-p^1$  satisfy  the following equations:
$$\left.\begin{array}{c} \pa_tv^2 -\De\,v^2=-\na\,p^2
    \\
  \div\,v^2=0
\end{array}\right\}\quad\mbox{in}\quad Q^+((1+\tau)/2)$$
$$
    v^2=0\qquad\mbox
    {on}\quad \Gamma((1+\tau)/2)\times ]-((1+\tau)/2)^2,0[.
$$
As it was shown in \ci{S1}, see Proposition 2 there,
\begin{equation}\label{310}
    v^2\in W^{2,1}_{s,\frac 32}(Q^+(\tau)),\qquad
    p^2\in W^{1,0}_{s,\frac 32}(Q^+(\tau))
\end{equation}
for any $s>9/8$.  (\ref{35}) follows from (\ref{39}),
(\ref{310}), and the obvious inequality $3/2>4/3$.

Since $v$ and $p$ are a suitable weak solution, we may apply
various conditions of the so-called $\varepsilon$-regularity.
First, we would like to note that  pairs $v$, $p$ and $v$,
$p-[p]_{B^+}$ are  suitable weak solutions to the Navier-Stokes
equations in $Q^+$ near $\Gamma\times [-1,0]$ simultaneously.
Therefore, the main result of \ci{S5}, see Theorem 1.2 in \ci{S5},
may be formulated in the following way.
\begin{lemma}\la{3l1}  There exist  universal positive constants
$\varepsilon_1$ and $c^1_k$, $k=0,1,...,$ with the following
property. Let a pair $v$ and $p$ be an arbitrary suitable weak
solution to the Navier-Stokes equations in $Q^+$ near
$\Gamma\times [-1,0]$ and satisfy the additional condition
\begin{equation}\label{311}
    C^+(0,1;v)+D^+(0,1;p)<\varepsilon_1.
\end{equation}
Then, for any $k=0,1,...$, the function $\na^k\,v$ is H\"older
continuous in the closure of the set ${Q}^+(1/2)$ and
$$\sup\limits_{z\in Q^+(1/2)}|\na^k\,v|\leq c^1_k.$$
 \end{lemma}

By the embedding theorem, we can reformulate Lemma \ref{3l1} in
the following way.
\begin{lemma}\la{3l2}  There exist  universal positive constants
$\varepsilon_2$ and $c^2_k$, $k=0,1,...,$ with the following
property. Let a pair $v$ and $p$ be an arbitrary suitable weak
solution to the Navier-Stokes equations in $Q^+$ near
$\Gamma\times [-1,0]$ and satisfy the additional condition
\begin{equation}\label{312}
    C^+(0,1;v)+D^+_1(0,1;p)<\varepsilon_2.
\end{equation}
Then, for any $k=0,1,2,...$, the function $\na^k\,v$ is H\"older
continuous in in the closure of the set ${Q}^+(1/2)$ and
$$\sup\limits_{z\in Q^+(1/2)}|\na^k\,v|\leq c^2_k.$$
 \end{lemma}

Finally, we would like to use another condition of
$\varepsilon$-regularity in terms
of the velocity $v$ only.
\begin{lemma}\la{3l3}  There exists an universal positive constant
$\varepsilon_3$ with the following property. Let a pair $v$ and
$p$ be an arbitrary suitable weak solution to the Navier-Stokes
equations in $Q^+$ near $\Gamma\times [-1,0]$. Assume that, for
some $R_0\in ]0,1[$, $v$ satisfies the additional condition
\begin{equation}\label{313}
   \sup\limits_{0<R\leq R_0}\frac
   1{R^2}\int\limits_{Q^+(R)}|v|^3\,dz<\varepsilon_3.
\end{equation}
Then, there exists $r_0\in ]0,R_0[$ such that  the function $v$ is
H\"older continuous in the closure of the set ${Q}^+(r_0)$.
 \end{lemma}
\textsc{Proof} Our arguments are similar to those used in the
proof of Proposition \ref{2p6}. By Lemma \ref{2l2}, we have
\begin{equation}\label{314}
D^+_1(r)\leq c\Big\{\Big(\frac
r\rho\Big)^2\Big[D^+_1(\rho)+(E^+)^\frac 34(\rho)\Big]+\Big(\frac
\rho r\Big)^3(A^+)^\frac 12(\rho)E^+(\rho)\Big\}
\end{equation}
for any $0<r\leq\rho\leq R_0/2$. Here,
\begin{eqnarray*}
  A^+(\rho)\equiv  A^+(0,\rho;v),\qquad  C^+(\rho)\equiv  C^+(0,\rho;v)\\
 D^+_1(\rho)\equiv  D^+_1(0,\rho;p),\qquad E^+(\rho)\equiv
 E^+(0,\rho;v).
\end{eqnarray*}
In addition, the local energy inequality gives us:
\begin{eqnarray}\la{315}
 \nonumber
 A^+(\rho)+ E^+(\rho)\leq c\Big[(C^+)^\frac 23(2\rho)+(C^+)^\frac 13(2\rho)
 (D^+)^\frac 23(2\rho)+C^+(2\rho)\Big]\\
   \leq c\Big[\varepsilon_3^\frac 23+\varepsilon_3^\frac 13(D^+_1)^\frac 23(2\rho)
   +\varepsilon_3\Big].
\end{eqnarray}
Without loss of generality, we may assume that $\varepsilon_3\leq
1$. Combining (\ref{314}) and (\ref{315}), we find
\begin{eqnarray*}
 D^+_1(r)\leq c\Big\{\Big(\frac
r\rho\Big)^2\Big[D^+_1(\rho)+\Big( \varepsilon_3^\frac
23+\varepsilon_3^\frac 13(D^+_1)^\frac 23(2\rho)
   +\varepsilon_3\Big )^\frac 34 \Big]\\
+\Big(\frac \rho r\Big)^3  \Big( \varepsilon_3^\frac
23+\varepsilon_3^\frac 13(D^+_1)^\frac 23(2\rho)
   +\varepsilon_3\Big )^\frac 32\Big\}   \\
 \leq c\Big\{\Big(\frac
r\rho\Big)^2\Big[D^+_1(2\rho) +\varepsilon_3^\frac 12\Big]+
\varepsilon_3^\frac 12  \Big(\frac \rho r\Big)^3D^+_1(2\rho)+
\varepsilon_3^\frac 12  \Big(\frac \rho r\Big)^3\Big\} \\
  \leq c\Big\{ \Big[\Big(\frac
r\rho\Big)^2+ \varepsilon_3^\frac 12  \Big(\frac \rho
r\Big)^3\Big]D^+_1(2\rho)+ \varepsilon_3^\frac 12  \Big(\frac \rho
r\Big)^3\Big\}
\end{eqnarray*}
for any   $0<r\leq\rho\leq R_0/2$. After some simple calculations,
the latter can be rewritten as follows.
$$D^+_1(\tau R) \leq c\Big\{ \Big[\tau^2+\frac{\varepsilon_3^\frac
12}{\tau^3}\Big]D^+_1(R)+\frac{\varepsilon_3^\frac
12}{\tau^3}\Big\}$$ for any $0<R\leq R_0$ and for any $0<\tau<1$.
Next, we fix $\tau$ so that
$$c\tau<1/2$$ and assume that
$$c\frac{\varepsilon_3^\frac
12}{\tau^3}<\frac \tau 2 \Leftrightarrow (\varepsilon_3<\Big(\frac
{\tau^4}{2c}\Big)^2).$$ Then we have
$$D^+_1(\tau R) \leq \tau D^+_1(R)+c\frac{\varepsilon_3^\frac
12}{\tau^3}$$ for any $0<R\leq R_0$. Making iterations, we find
$$D^+_1(\tau^k R_0) \leq \tau^k D^+_1(R_0)+c\frac{\varepsilon_3^\frac
12}{\tau^3}\,\frac 1{1-\tau}$$ for any natural $k$. Therefore,
\begin{eqnarray*}
 C^+(\tau^k R_0)+D^+_1(\tau^k R_0)\leq\varepsilon_3+\tau^k
D^+_1(R_0)+c\frac{\varepsilon_3^\frac 12}{\tau^3}\,\frac
1{1-\tau}  \\
\leq \tau^k D^+_1(R_0)+2c\frac{\varepsilon_3^\frac
12}{\tau^3}\,\frac 1{1-\tau}
\end{eqnarray*}
Choose $\varepsilon_3$ so small that
$$2c\frac{\varepsilon_3^\frac
12}{\tau^3}\,\frac 1{1-\tau}<\frac {\varepsilon_2}3$$ and then fix
$k$ so  that
$$\tau^k
D^+_1(R_0)<\frac {\varepsilon_2}3.$$ Hence, a
$\varepsilon$-regularity condition holds. In particular, we have
$$C^+(\tau^k R_0)+D^+_1(\tau^k R_0)<\varepsilon_2.$$
By scaling and Lemma \ref{3l2}, we can take $r_0=\frac 12 \tau^k
R_0$. Lemma \ref{3l3} is proved.

\setcounter{equation}{0}
\section{Proof of Theorem \ref{it1} }

We let
\begin{equation}\label{41}
    L\equiv \|v\|_{3,\infty,Q^+}<+\infty.
\end{equation}
Using known arguments, (\ref{35}), and (\ref{41}), we can assert
that
\begin{equation}\label{42}
    \sup\limits_{-(3/4)^2\leq t\leq
    0}\|v(\cdot,t)\|_{3,B^+(3/4)}\leq L.
\end{equation}

Assume now that the statement of Theorem \ref{it1} is false. Let
$z_0\in \overline{Q}^+(1/2)$ be a singular point. As it was shown
in \ci{ESS4}, Theorem 1.4, $z_0$ must belong to
$\overline{\Gamma}(1/2)$. Without loss of generality (just by
translation and by scaling), we may assume that $z_0=0$. It
follows from Lemma \ref{3l3} that a sequence $R_k\downarrow 0$
exists such that
\begin{equation}\label{43}
    \frac 1{R_k}\int\limits_{Q^+(R_k)}|v|^3\,dz\geq\varepsilon_3
\end{equation}
for any natural $k$.

Extending functions $v$ and $p$ outside $Q^+$ to zero, we
introduce scaled functions
$$u^k(y,s)=R_kv(R_ky,R^2_ks),\qquad q^k(y,s)=R^2_kp(R_ky,R^2_ks)$$
for $y\in \mathbb R^3_+$ and for $s\in \mathbb R_-=\{s<0\}$.

Our first observation is that
\begin{equation}\label{44}
    u^k{\stackrel{\star}{\rightharpoonup}}\,u\qquad\mbox {in}\qquad
L_{3,\infty}(\mathbb R^3_+\times \mathbb R_-)
\end{equation}
(at least for a subsequence).

Fix $a>0$ and let $k$ be so that
\begin{equation}\label{45}
    R_ka<\frac 18.
\end{equation}

By Proposition 1 in \ci{S1}, we have two estimates:
\begin{eqnarray}\la{46}
  \nonumber  \|\na^2\,u^k\|_{\frac 98,\frac 32, Q^+(a)}
  + \|\na\,q^k\|_{\frac 98,\frac 32, Q^+(a)}\\
\nonumber
\leq c_1(a) \Big[\|u^k_iu^k_{,i}\|_{\frac 98,\frac 32,
Q^+(2a)}+
\|u^k\|_{W^{1,0}_{\frac 98,\frac 32}( Q^+(2a))}  \\
\nonumber
 + \|q^k-[q^k]_{B^+(2a)}\|_{\frac 98,\frac 32, Q^+(2a)}\Big]\\
 \nonumber
  \leq  c^\prime_1(a) \Big[\|u^k\|_{3,\infty, Q^+(2a)} \|\na\,
  u^k\|_{2, Q^+(2a)}+\|\na\,
  u^k\|_{2, Q^+(2a)}\\
  \nonumber
  +\|u^k\|_{3, Q^+(2a)}+
   \|q^k-[q^k]_{B^+(2a)}\|_{\frac 98,\frac 32, Q^+(2a)}\Big]
 \\
 \nonumber
 \leq  c^{\prime\prime}_1(a,L) \Big[\|\na\,u^k\|^2_{2, Q^+(2a)}+1
  +
   \|q^k-[q^k]_{B^+(2a)}\|^\frac 32_{\frac 98,\frac 32, Q^+(2a)}\Big] \\
  \leq  c^{\prime\prime\prime}_1(a,L) \Big[ E^+(0,2a;u^k)+1+D_1^+(0,2a;q^k)
   \Big]
\end{eqnarray}
and
\begin{eqnarray}\la{47}
  \nonumber
   \|\na^2\,u^k\|_{\frac 43, Q^+(a)}
  + \|\na\,q^k\|_{\frac 43, Q^+(a)}\\
  \nonumber
  \leq c_2(a) \Big[\|u^k_iu^k_{,i}\|_{\frac 43,
Q^+(2a)}+
\|u^k\|_{W^{1,0}_{\frac 43}( Q^+(2a))}  \\
\nonumber
 + \|q^k-[q^k]_{B^+(2a)}\|_{\frac 43, Q^+(2a)}\Big]\\
  \nonumber
  \leq  c^\prime_2(a,L) \Big[ \|\na\,
  u^k\|^\frac 32_{2, Q^+(2a)}+\|\na\,
  u^k\|_{2, Q^+(2a)}
  \\
  \nonumber
  +1+
   \|q^k-[q^k]_{B^+(2a)}\|_{\frac 32, Q^+(2a)}\Big]\\
 \leq  c^{\prime\prime}_2(a,L) \Big[ E^+(0,2a;u^k)+1+D_1^+(0,2a;q^k)
   \Big].
\end{eqnarray}
On the other hand, by the inverse scaling and by the local energy
inequality, we find (see the proof of Proposition \ref{2p6})
\begin{eqnarray}\la{48}
  \nonumber
   E^+(0,2a;u^k)+D_1^+(0,2a;q^k) =
    E^+(0,2aR_k;v)+D_1^+(0,2aR_k;p)\\
\leq c_3(L)\Big[1+D_1^+(0,4aR_k;p)\Big].
\end{eqnarray}
 To establish uniform
boundedness with respect to $k$, let us make use of  Proposition
\ref{2p6}. As a result, we have
\begin{eqnarray}\la{49}
  \nonumber
  D_1^+(0,4aR_k;p)\leq c_4(L)\Big[(4aR_k)D_1^+(0,1;p)+1\Big]
  \\
\leq c_4(L)\Big[D_1^+(0,1;p)+1\Big].
\end{eqnarray}
Now, by (\ref{46})--(\ref{49}) and by the diagonal Cantor process,
we can select subsequences (still denoted by $u^k$ and $q^k$) with
the following properties:
\begin{eqnarray}
\la{410}  u^k\rightharpoonup\,u\qquad\mbox{in}\quad W^{2,1}_{\frac
43}(Q^+(a))\cap
 W^{2,1}_{\frac 98,\frac 32}(Q^+(a)),  \\
  \la{411} q^k\rightharpoonup\,q\qquad\mbox{in}\quad W^{1,0}_{\frac 43}(Q^+(a))\cap
 W^{1,0}_{\frac 98,\frac 32}(Q^+(a)),  \\
\la{412} \na\,u^k\rightharpoonup\,\na\,u\qquad\mbox{in}\quad
L_2(Q^+(a))
\end{eqnarray}
for any $a>0$. Moreover, by (\ref{44}), (\ref{410}), (\ref{411}),
and by the known multiplicative inequality, we can state that:
\be\la{413}\left.\begin{array}{c}u^k\rightharpoonup\,u
\qquad\mbox{in}\quad L_\frac {10}3(Q^+(a)
)\\
u^k\rightarrow\,u \qquad\mbox{in}\quad L_3(Q^+(a)
)\end{array}\right\}\ee for any $a>0$.

According to (\ref{410})--(\ref{413}), the pair $u$ and $q$ forms
a suitable weak solution to the Navier-Stokes equations in
$Q^+(a)$ near $\Gamma(a)\times [-a^2,0]$ for any $a>0$. This
solution possesses the additional property
\begin{equation}\label{414}
    \|u\|_{3,\infty,\mathbb R^3_+\times \mathbb R_-}\leq  L.
\end{equation}
Moreover, using (\ref{410}) and interpolation, we can show that
\begin{equation}\label{415}
u^k\rightarrow\,u \qquad\mbox{in}\quad C([-a^2,0];L_2(B^+(a))),
\end{equation}
see details in the proof of (3.23) in \ci{ESS4}. Letting
$d=x_{03}/2$ for an arbitrary point $x_0\in \mathbb R^3_+$ and
using (\ref{42}) and (\ref{415}), we find
\begin{eqnarray*}
  \Big(\int\limits_{B(x_0,d)}|u(x,0)|^2\,dx\Big)^\frac 12
  \leq \Big(\int\limits_{B(x_0,d)}|u(x,0)-u^k(x,0)|^2\,dx\Big)^\frac 12 \\
  + \Big(\int\limits_{B(x_0,d)}|u^k(x,0)|^2\,dx\Big)^\frac 12
  \leq \Big(\int\limits_{B(x_0,d)}|u(x,0)-u^k(x,0)|^2\,dx\Big)^\frac 12 \\
+cd\Big(\int\limits_{B(x_0,d)}|u^k(x,0)|^3\,dx\Big)^\frac 23
\\
  = \Big(\int\limits_{B(x_0,d)}|u(x,0)-u^k(x,0)|^2\,dx\Big)^\frac 12
  +cd\Big(\int\limits_{B(x_0R_k,dR_k)}|v(x,0)|^3\,dx\Big)^\frac 23\rightarrow
  \,0
\end{eqnarray*}
as $k\to +\infty$. So,
\begin{equation}\label{416}
    u(\cdot,0)=0\qquad\mbox{in}\quad \mathbb R^3_+.
\end{equation}
However, $u$ is not trivial solution. This directly follows from
(\ref{43}) and (\ref{413}):
\begin{equation}\label{417}
    \frac 1{R^2_k}\int\limits_{Q^+(R_k)}|v|^3\,dz=
    \int\limits_{Q^+}|u^k|^3\,dz\to\,\int\limits_{Q^+}|u|^3\,dz\geq
    \varepsilon_3.
\end{equation}

Now, our goal is to show that in fact
\begin{equation}\label{418}
    u=0\qquad\mbox{in}\quad\mathbb R^3_+\times ]-1,0[.
\end{equation}
This would contradict with (\ref{417}) and  complete the proof of
Theorem \ref{it1}.

First, we would like to point out that by interior regularity
results of \ci{ESS4}, see Theorem 1.4 there, $u$ may have singular
points on the plane $x_3=0$ only. To apply backward uniqueness
arguments (as it was done in \ci{ESS4}), we need to know if
$u$ and $\na\,u$ are bounded on certain sets. We shall show that
it is so on sets of the form $(\mathbb R^3_++hi_3)\times ]-T,0[$,
where $h>0$ and $T>0$ are arbitrarily fixed and $i_3=(0,0,1)$. To
this end, let us prove the following statement.
\begin{lemma}\la{4l1}
There exists a positive constant $c_5$, depending on  $L$ and
$D^+_1(0,\\1;p)$ only, with the following property.
Fix $h>0$ and $T>0$ arbitrarily, then
\begin{equation}\label{419}
    D(e_0,2h;q)\leq c_5
\end{equation}
for any $e_0=(y_0,s_0)\in (\mathbb R^3_++3hi_3)\times ]-T,0[$.
\end{lemma}
\textsc{Proof} From (\ref{411}), we know that
\begin{equation}\label{420}
    \limsup\limits_{k\to +\infty} D(e_0,2h;q^k)\geq D(e_0,2h;q).
\end{equation}
So, it is sufficient to prove the following bound
\begin{equation}\label{421}
D(e_0,2h;q^k)\leq c_5(L,D^+_1(0,1;p))
\end{equation}
provided that
\begin{equation}\label{422}
    x_0^k=y_0R_k\in B^+(1/4),\qquad t_0^k=s_0R^2_k>-(1/4)^2.
\end{equation}
Obviously,  (\ref{420}) and (\ref{421}) imply (\ref{419}).

We have
\begin{equation}\label{423}
    D(e_0,2h;q^k)=D(z^k_0,2hR_k;p),\qquad z_0^k=( x_0^k,t_0^k).
\end{equation}
Further,  if $d_k=x^k_{03}=\dist(x_0^k,\Gamma)>2hR_k$, then
$Q(z_0^k,d_k)\subset Q^+(1/2)$ and, therefore, we may use
Proposition \ref{2p5}. As a result, we find
\begin{eqnarray}\la{424}
  \nonumber
 D(z^k_0,2hR_k;p)\leq c_6(L)\Big[\Big(\frac {2hR_k}{d_k}\Big)^\frac 54
   D(z^k_0,d_k;p)+1\Big]\\
\leq c_6(L)\Big[
   D(z^k_0,d_k;p)+1\Big].
\end{eqnarray}
On the other hand,
$$Q(z_0^k,d_k)\subset Q^+(\overline{z}_0^k,2d_k),\qquad
\overline{z}_0^k=( \overline{x}_0^k,t_0^k)$$ where
$\overline{x}_0^k=(x_{01}^k,x_{02}^k,0)$, and, moreover,
$$ D(z^k_0,d_k;p)\leq cD^+(\overline{z}_0^k,2d_k;p).$$
Therefore, we have (see (\ref{424}))
\begin{eqnarray}\la{425}
  \nonumber
  D(z^k_0,2hR_k;p)\leq
  c_7(L)\Big[D^+(\overline{z}_0^k,2d_k;p)+1\Big]
  \\
\leq  c^\prime_7(L)\Big[D_1^+(\overline{z}_0^k,2d_k;p)+1\Big].
\end{eqnarray}
Now, taking into account
$$Q^+(\overline{z}_0^k,2d_k)\subset
Q^+(\overline{z}_0^k,1/2)\subset Q^+,$$ we apply Proposition
\ref{2p6} which says that
\begin{eqnarray}\la{426}
  \nonumber
D_1^+(\overline{z}_0^k,2d_k;p)\leq
c_8(L)\Big[\Big(\frac{2d_k}{1/2}\Big)
 D_1^+(\overline{z}_0^k,1/2;p)+1\Big] \\
\leq c_8(L)\Big[
 D_1^+(\overline{z}_0^k,1/2;p)+1\Big]
\leq c^\prime_8(L)\Big[
 D_1^+(0,1;p)+1\Big].
\end{eqnarray}
Obviously, (\ref{421}) follows from (\ref{423}), (\ref{425}), and
(\ref{426}). Lemma \ref{4l1} is proved.

Now, we proceed the proof of Theorem \ref{it1}. Fix $h\in ]0,1[$
arbitrarily and  let $T=100$. Take an arbitrary point
$z_0=(x_0,t_0)$ so that
$$z_0\in (\mathbb R^3_++3hi_3)\times ]-100,0[.$$
In the ball $B(x_0,2h)$, we decompose pressure
$$q=q_1+q_2$$
in the following way:
$$\int\limits_{B(x_0,2h)}q_1(x,t)\De\varphi(x)\,dx=
-\int\limits_{B(x_0,2h)}u(x,t)\otimes u(x,t):\na^2\varphi(x)\,dx$$
for any $\varphi\in C^2(\overline{B}(x_0,2h))$ such that
$\varphi|_{\pa B(x_0,2h)}=0$, and
$$\De q_2(\cdot,t)=0\qquad\mbox{in}\quad B(x_0,2h).$$
For $q_1$ and $q_2$, the following estimates are valid:
\begin{equation}\label{427}
\int\limits_{B(x_0,2h)}|q_1(x,t)|^\frac 32\,dx \leq c
\int\limits_{B(x_0,2h)}|u(x,t)|^3\,dx
\end{equation}
and
\begin{equation}\label{428}
    \sup\limits_{x\in B(x_0,h)}|\na\,q_2(x,t)|\leq c
   \frac 1h \Big(\frac 1{h^3}\int\limits_{B(x_0,2h)}|q_2(x,t)-
    [q_2]_{B(x_0,2h)}(t)|^\frac 32\,dx\Big)^\frac 23.
\end{equation}

For any $0<\rho <1$, using (\ref{427}) and (\ref{428}), we can
derive the  estimate:
\begin{eqnarray*}
   D(z_0,h\rho;q)\leq c\Big[  D(z_0,h\rho;q_1)+D(z_0,h\rho;q_2)\Big]\\
   \leq c\Big[\frac 1{(h\rho)^2}\int\limits_{Q(z_0,2h)}|q_1|^\frac 32\,dz
   +(h\rho)^{\frac 52}\int\limits_{t_0-(h\rho)^2}^{t_0}
   \sup_{x\in B(x_0,h)}| \na\,q_2(x,t)|^\frac 32\,dt\Big] \\
   \leq c\Big[\frac 1{(h\rho)^2}\int\limits_{Q(z_0,2h)}|u|^ 3\,dz
   +\rho^\frac 52 D(z_0,2h;q_2)\Big] \\
  \leq c\Big[\frac 1{(h\rho)^2}\int\limits_{Q(z_0,2h)}|u|^ 3\,dz
   +\rho^\frac 52 D(z_0,2h;q)+\rho^\frac 52 D(z_0,2h;q_1)\Big]  \\
  \leq c\Big[\frac 1{(h\rho)^2}\int\limits_{Q(z_0,2h)}|u|^ 3\,dz
   +\rho^\frac 52 D(z_0,2h;q)\Big].
\end{eqnarray*}
To evaluate the last term on the right hand side of the latter
inequality, we need Lemma \ref{4l1}. It gives us:
\begin{eqnarray*}
  C(z_0,h\rho;v)+D(z_0,h\rho;q)\leq c
 \Big[\frac 1{(h\rho)^2}\int\limits_{Q(z_0,2h)}|u|^ 3\,dz
   +\rho^\frac 52c_5\Big]
\end{eqnarray*}
for any  $z_0\in (\mathbb R^3_++3hi_3)\times ]-100,0[$.

Now,  fix
$\rho(L,D^+_1(0,1;p))\in ]0,1[$ in the following way:
$$c\rho^\frac 52c_5<\varepsilon_3/3.$$
Then, we find $R_1>100
$ so  that
$$\frac 1{(h\rho)^2}\int\limits_{-200}^0\,dt\int\limits_{\mathbb
R^3_+\setminus B(R_1/4)}|u|^3\,dx<\varepsilon_3/3.$$ Since
$Q(z_0,2h)\subset(\mathbb R^3_+\setminus B(R_1/4))\times ]-200,0[$
for $|x_0|>R_1/2$, the latter implies
$$\frac 1{(h\rho)^2}\int\limits_{Q(z_0,2h)}|u|^ 3\,dz<\varepsilon_3/3$$
for all $z_0\in (\mathbb R^3_++3hi_3)\times ]-100,0[$ such that
$|x_0|>R_1/2$. It is known (see, for instance, \ci{ESS4}, Lemma
2.2) that, for any $k=0,1,...$, the function $\na^k\,u$ is bounded
on the set $(\mathbb R^3_++6hi_3)\setminus B(R_1)\times [-50,0]$.
Smoothness (and boundedness) of $\na^k\,u$ on the set $(\mathbb
R^3_++6hi_3)\cap\overline{B}(R_1)\times [-50,0]$ is already known.

So, if we introduce the vorticity $\omega=\na\wedge u$, then
$\omega$ satisfies the  relations:
\begin{eqnarray*}
  |\pa_t\omega-\De\omega|\leq  M(|\omega|+|\na\omega|),\\
  |\omega|\leq M
\end{eqnarray*}
on the set $(\mathbb R^3_++6hi_3)\times [-50,0]$ for some $M>0$,
and
$$\omega(\cdot,0)=0\qquad \mbox{in}\quad \mathbb R^3_+.$$
In \ci{ESS3}, it was shown that these three conditions imply
$$\omega=0\qquad \mbox{in}\quad (\mathbb R^3_++6hi_3)\times
[-50,0].$$ Since $h$ was taken arbitrarily, the latter means that
$$\omega=0\qquad \mbox{in}\quad \mathbb R^3_+\times
[-50,0].$$ Hence, for a.a. $t\in [-50,0]$, $u$ is a harmonic
function, which satisfies the boundary condition  $u(x,t)=0$ if
$x_3=0$. But, for a.a. $t\in [-50,0]$, $L_3$-norm of $u$ over
$\mathbb R^3_+$ is finite.  This leads  to the  conclusion that,
for the same $t$, $u(\cdot,t)=0$ in $\mathbb R^3_+$ . Theorem
\ref{it1} is proved.

\setcounter{equation}{0}
\section{Application to the initial boundary value problem in a half space }

Fix an arbitrary $T>0$ and consider the following initial boundary
value problem: \be\la{51}\left.\begin{array}{c}
\pa_tv+\div\,v\otimes v-\De\,v=-\na\,p
    \\
  \div\,v=0
\end{array}\right\}\quad\mbox{in}\quad Q_T=\mathbb R^3_+\times ]0,T[;\ee
\begin{equation}\label{52}
    v(x,t)=0,\qquad x_3=0\quad\mbox{and}\quad 0\leq t\leq T;
\end{equation}
\begin{equation}\label{53}
    v(x,0)=a\qquad x\in \mathbb R^3_+,
\end{equation}
where a solenoidal vector-valued field $a$ belongs to $L_2(\mathbb
R^3_+)$. For any $T>0$,  problem (\ref{51})--(\ref{53}) has at
least one the so-called weak Leray-Hopf solution $v$ having the
following properties (see, for instance, \ci{L1} and  \ci{L}):
$$v\in L_{2,\infty}(Q_T)\cap W^{1,0}_2(Q_T);$$
$$t\mapsto \int\limits_{\mathbb R^3_+}v(x,t)\cdot
u(x)\,dx\,\,\mbox{is continuous on}\,\,[0,T]\,\,\mbox{for any}\,\,
u\in L_2(\mathbb R^3_+);$$
$$\int\limits_{Q_T}(v\cdot \pa_tw+v\otimes
v:\na\,w-\na\,v:\na\,w)\,dz=0\,\,$$ for any divergence free
 test function $w\in C^\infty_0(Q_T)$;
 $$\|v(\cdot,t)-a(\cdot)\|_{L_2(\mathbb R^3_+)}\to\,0$$
 as $t\downarrow 0$;
 $$\int\limits_{\mathbb R^3_+}|v(x,t)|^2\,dx+2\int\limits_{\mathbb
 R^3_+\times ]0,t[}|\na\,v|^2\,dxdt^\prime\leq\int
 \limits_{\mathbb R^3_+}|a(x)|^2\,dx$$
 for any $t\in [0,T]$.
 \begin{theorem}\la{5t1} Assume that
$$
    v\in L_{3,\infty}(Q_T).
$$
Then, $v\in L_5(Q_T)$ and, moreover, $v$ is smooth and unique.
\end{theorem}
\textsc{Proof} To this end, as it was shown in \ci{ESS4},  it
sufficient to prove the  estimate
$$\sup\limits_{z\in\mathbb R^3_+\times [\de,T]}|v(z)|\leq M(\de)<
\infty,\qquad\forall\delta>0.$$  With the help of  linear theory,
the associated pressure $p$ can be introduced so that:
$$p\in L_\frac 32(B^+(R)\times ]\de,T[)$$
for any $R>0$;$$ \na^2\,v, \,\pa_t\,v,\,\na\,p\in L_{\frac
98,\frac 32}(\mathbb R^3_+\times ]\de,T[);$$
$$\pa_tv+\div\,v\otimes v-\De\,v=-\na\,p$$
a.a. in $Q_T$.

From Theorem \ref{it1} and from Theorem 1.4 in \ci{ESS4}, it
follows that $v$ has no singular point. We must prove the global
boundedness only.  Obviously, for $R\to\,+\infty$,
$$\int\limits_0^T\,dt\int\limits_{\mathbb R^3_+\setminus
B^+(R)}|v|^3\,dx\to\,0$$ and
$$\int\limits_\de^T\,dt\Big(\int\limits_{\mathbb R^3_+\setminus
B^+(R)}|\na\,p|^\frac 98\,dx\Big)^\frac 43\to\,0.$$

Next,  using these facts as well as various conditions of
$\varepsilon$-regularity, see Lemmas \ref{3l1} and \ref{3l2} and
Lemma 2.2 in \ci{ESS4}, we conclude that
$$\sup\limits_{z\in(\mathbb R^3_+\setminus B^+(R))\times
[\de,T]}|v(z)|\leq M_1(\de,R)<\infty.$$ Another estimate
$$\sup\limits_{z\in\overline{ B}^+(R)\times [\de,T]}
|v(z)|\leq M_2(\de,R)<\infty$$ is already known since our solution
$v$ is locally smooth. Theorem \ref {5t1} is proved.

\noindent
G. Seregin\\
Steklov Institute of Mathematics at St.Petersburg, \\
St.Petersburg, Russia

\end{document}